\documentclass[12pt]{amsart}
\usepackage{amssymb}
\title[Irregular manifolds]{Irregular manifolds with a canonical
linear system, composite with a pencil}
\author{Jin-Xing Cai and Eckart Viehweg}
\address{School of  Mathematical Sciences  and Institute of Mathematics,
Peking University, Beijing 100871, P. R. China}
\email{cai@math.pku.edu.cn}
\address{Universit\"at Essen, FB6 Mathematik, 45117 Essen, Germany}
\email{viehweg@uni-essen.de}
\thanks{This work has been supported by the ``DFG-Schwerpunktprogramm
Globale Methoden in der Komplexen Geometrie'' and by the DFG-NSFC
Chinese-German project ``Komplexe Geometrie''. The first named
author is also partially supported by
the NSFC (No. 10271005) and SRF for ROCS, SEM}
\newcommand{\op}{\mathrm}
\newcommand{\Cal}{\mathcal}
\newcommand{\mc}{\mathcal}
\newcommand{\ti}{\tilde}
\newcommand{\BB}{\mathbb}

\theoremstyle{plain}
\newtheorem{theorem}{Theorem}[section]
\newtheorem{thm}[theorem]{Theorem}
\newtheorem{prop}[theorem]{Proposition}
\newtheorem{lem}[theorem]{Lemma}
\newtheorem{cor}[theorem]{Corollary}
\theoremstyle{definition}
\newtheorem{rk}[theorem]{Remark}

\newtheorem{exa}[theorem]{Example}

\catcode`\@=11
\def\opn#1#2{\def#1{\mathop{\kern0pt\fam0#2}\nolimits}}
\def\underrightarrow{\mathpalette\underrightarrow@}
\def\underrightarrow@#1#2{\vtop{\ialign{$##$\cr
 \hfil#1#2\hfil\cr\noalign{\nointerlineskip}%
 #1{-}\mkern-6mu\cleaders\hbox{$#1\mkern-2mu{-}\mkern-2mu$}\hfill
 \mkern-6mu{\to}\cr}}}

\def\underleftarrow{\mathpalette\underleftarrow@}
\def\underleftarrow@#1#2{\vtop{\ialign{$##$\cr
 \hfil#1#2\hfil\cr\noalign{\nointerlineskip}#1{\leftarrow}\mkern-6mu
 \cleaders\hbox{$#1\mkern-2mu{-}\mkern-2mu$}\hfill
 \mkern-6mu{-}\cr}}}
\let\amp@rs@nd@\relax
\newdimen\ex@
\newdimen\bigaw@
\newdimen\minaw@
\newdimen\minCDaw@
\newif\ifCD@
\def\minCDarrowwidth#1{\minCDaw@#1}
\newenvironment{CD}{\@CD}{\@endCD}
\def\@CD{\def\A##1A##2A{\llap{$\vcenter{\hbox
 {$\scriptstyle##1$}}$}\Big\uparrow\rlap{$\vcenter{\hbox{%
$\scriptstyle##2$}}$}&&}%
\def\V##1V##2V{\llap{$\vcenter{\hbox
 {$\scriptstyle##1$}}$}\Big\downarrow\rlap{$\vcenter{\hbox{%
$\scriptstyle##2$}}$}&&}%
\def\={&\hskip.5em\mathrel
 {\vbox{\hrule width\minCDaw@\vskip3\ex@\hrule width
 \minCDaw@}}\hskip.5em&}%
\def\verteq{\Big\Vert&&}%
\def\noarr{&&}%
\def\vspace##1{\noalign{\vskip##1\relax}}\relax\let\amp@rs@nd@&\iffalse}\fi
 \CD@true\vcenter\bgroup\relax\let\\=\cr\iffalse}\fi\tabskip\z@skip\baselineskip20\ex@
 &\hfill$\m@th##$\hfill\cr}
\def\@endCD{\cr\egroup\egroup}
\def\>#1>#2>{\amp@rs@nd@\setbox\z@\hbox{$\scriptstyle
 \;{#1}\;\;$}\setbox\@ne\hbox{$\scriptstyle\;{#2}\;\;$}\setbox\tw@
 \hbox{$#2$}\ifCD@
 \global\bigaw@\minCDaw@\else\global\bigaw@\minaw@\fi
 \ifdim\wd\z@>\bigaw@\global\bigaw@\wd\z@\fi
 \ifdim\wd\@ne>\bigaw@\global\bigaw@\wd\@ne\fi
 \ifCD@\hskip.5em\fi
 \ifdim\wd\tw@>\z@
 \mathrel{\mathop{\hbox to\bigaw@{\rightarrowfill}}\limits^{#1}_{#2}}\else
 \mathrel{\mathop{\hbox to\bigaw@{\rightarrowfill}}\limits^{#1}}\fi
 \ifCD@\hskip.5em\fi\amp@rs@nd@}
\def\<#1<#2<{\amp@rs@nd@\setbox\z@\hbox{$\scriptstyle
 \;\;{#1}\;$}\setbox\@ne\hbox{$\scriptstyle\;\;{#2}\;$}\setbox\tw@
 \hbox{$#2$}\ifCD@
 \global\bigaw@\minCDaw@\else\global\bigaw@\minaw@\fi
 \ifdim\wd\z@>\bigaw@\global\bigaw@\wd\z@\fi
 \ifdim\wd\@ne>\bigaw@\global\bigaw@\wd\@ne\fi
 \ifCD@\hskip.5em\fi
 \ifdim\wd\tw@>\z@
 \mathrel{\mathop{\hbox to\bigaw@{\leftarrowfill}}\limits^{#1}_{#2}}\else
 \mathrel{\mathop{\hbox to\bigaw@{\leftarrowfill}}\limits^{#1}}\fi
 \ifCD@\hskip.5em\fi\amp@rs@nd@}

\def\@CDS{\def\A##1A##2A{\llap{$\vcenter{\hbox
 {$\scriptstyle##1$}}$}\Big\uparrow\rlap{$\vcenter{\hbox{%
$\scriptstyle##2$}}$}&}%
\def\V##1V##2V{\llap{$\vcenter{\hbox
 {$\scriptstyle##1$}}$}\Big\downarrow\rlap{$\vcenter{\hbox{%
$\scriptstyle##2$}}$}&}%
\def\={&\hskip.5em\mathrel
 {\vbox{\hrule width\minCDaw@\vskip3\ex@\hrule width
 \minCDaw@}}\hskip.5em&}
\def\verteq{\Big\Vert&}
\def\novarr{&}
\def\noharr{&&}
\def\SE##1E##2E{\slantedarrow(0,18)(4,-3){##1}{##2}&}
\def\SW##1W##2W{\slantedarrow(24,18)(-4,-3){##1}{##2}&}
\def\NE##1E##2E{\slantedarrow(0,0)(4,3){##1}{##2}&}
\def\NW##1W##2W{\slantedarrow(24,0)(-4,3){##1}{##2}&}
\def\slantedarrow(##1)(##2)##3##4{%
\thinlines\unitlength1pt\lower 6.5pt\hbox{\begin{picture}(24,18)%
\put(##1){\vector(##2){24}}%
\put(0,8){$\scriptstyle##3$}%
\put(20,8){$\scriptstyle##4$}%
\end{picture}}}
\def\vspace##1{\noalign{\vskip##1\relax}}\relax\let\amp@rs@nd@&\iffalse}\fi
 \CD@true\vcenter\bgroup\relax\let\\=\cr\iffalse}\fi\tabskip\z@skip\baselineskip20\ex@
 &\hfill$\m@th##$\hfill\cr}
\def\@endCDS{\cr\egroup\egroup}
\newdimen\TriCDarrw@
\newif\ifTriV@

\def\@TriCDV{\TriV@true\def\TriCDpos@{6}\@TriCD}
\def\@TriCDA{\TriV@false\def\TriCDpos@{10}\@TriCD}
\def\@TriCD#1#2#3#4#5#6{%
\setbox0\hbox{$\ifTriV@#6\else#1\fi$}
\TriCDarrw@=\wd0 \advance\TriCDarrw@ 24pt
\advance\TriCDarrw@ -1em
\def\SE##1E##2E{\slantedarrow(0,18)(2,-3){##1}{##2}&}
\def\SW##1W##2W{\slantedarrow(12,18)(-2,-3){##1}{##2}&}
\def\NE##1E##2E{\slantedarrow(0,0)(2,3){##1}{##2}&}
\def\NW##1W##2W{\slantedarrow(12,0)(-2,3){##1}{##2}&}
\def\slantedarrow(##1)(##2)##3##4{\thinlines\unitlength1pt
\lower 6.5pt\hbox{\begin{picture}(12,18)%
\put(##1){\vector(##2){12}}%
\put(-4,\TriCDpos@){$\scriptstyle##3$}%
\put(12,\TriCDpos@){$\scriptstyle##4$}%
\end{picture}}}
\def\={\mathrel {\vbox{\hrule
   width\TriCDarrw@\vskip3\ex@\hrule width
   \TriCDarrw@}}}
\def\>##1>>{\setbox\z@\hbox{$\scriptstyle
 \;{##1}\;\;$}\global\bigaw@\TriCDarrw@
 \ifdim\wd\z@>\bigaw@\global\bigaw@\wd\z@\fi
 \hskip.5em
 \mathrel{\mathop{\hbox to \TriCDarrw@
{\rightarrowfill}}\limits^{##1}}
 \hskip.5em}
\def\<##1<<{\setbox\z@\hbox{$\scriptstyle
 \;{##1}\;\;$}\global\bigaw@\TriCDarrw@
 \ifdim\wd\z@>\bigaw@\global\bigaw@\wd\z@\fi
 \mathrel{\mathop{\hbox to\bigaw@{\leftarrowfill}}\limits^{##1}}
 }
 \CD@true\vcenter\bgroup\relax\let\\=\cr\iffalse}\fi
 \tabskip\z@skip\baselineskip20\ex@
 \ifTriV@
 \halign\bgroup
 &\hfill$\m@th##$\hfill\cr
#1&\multispan3\hfill$#2$\hfill&#3\\
&#4&#5\\
&&#6\cr\egroup%
\else
 \halign\bgroup
 &\hfill$\m@th##$\hfill\cr
&&#1\\%
&#2&#3\\
#4&\multispan3\hfill$#5$\hfill&#6\cr\egroup
\fi}
\def\@endTriCD{\egroup}
\begin{document}
\begin{abstract}
Let $X$ be a complex projective $n$-dimensional manifold of
general type whose canonical system is composite with a pencil.
If the Albanese map is generically finite, but not surjective, or
if the irregularity is strictly larger than $n$ and the image of
$X$ in ${\rm Alb}(X)$ of Kodaira dimension one, then the geometric
genus $p_g(F)$ of a general fibre $F$ of the canonical map is one
and the latter factors through the Albanese map. The last part of this result
holds true for any threefold with $q(X)\geq 5$.
\end{abstract}
\maketitle
\section*{Introduction}\label{sec1}
In this note we study irregular higher dimensional complex
projective manifolds $X$ whose canonical map $\Phi_X:X \to
\BB{P}^{p_g(X)-1}$ has a one dimensional image. Hence replacing
$X$ by some blowing up, the canonical map factors through a
surjection $f: X \to C$ to a curve $C$ with connected fibres. We
will call $f:X\to C$ the canonical fibration of $X$.

As shown by  G. Xiao  \cite{X1} if the canonical system of a
surface $S$ of general type is composite with a pencil, then
$q(S)\leq 2$. Moreover this pencil is rational if and only if
$q(S)=2$.

Before formulating some generalizations to higher dimensional manifolds,
let us recall that the Albanese dimension  $\op{a} (X)$
of a manifold $X$, is defined to be the dimension of the image of
the Albanese map $\op{alb}=\op{alb}_X : X\to \op{Alb}(X)$.
$X$ is said to be of maximal Albanese dimension if $ a(X)=\dim(X)$.
Remark that the Kodaira dimension $\kappa(\op{alb}(X))$ is always
supposed to be the Kodaira dimension of a nonsingular model of $\op{alb}(X)$.

\begin{thm}[Section \ref{sect2} and Corollary \ref{cor}]\label{thm1}
Let $X$ be a projective  manifold of maximal Albanese
dimension. Assume the canonical system
is composite with a pencil, and let $f: X\to C$ be the canonical fibration,
with general fibre $F$.
\begin{enumerate}
\item[(i)] If $q(X)>\dim(X)$, then
\begin{eqnarray*}
& & p_g(X)=g(C) \ge 2, \ \ q(X)=g(C)+\dim(X)-1, \text{ \ and \ }
\chi(\mathcal O_X)=0.\\
& & p_g(F)=1, \text{ \ and \ } \chi(\mc O_F)=0.
\end{eqnarray*}
\item[(ii)] If $q(X)=\dim(X)$, then $C\simeq \mathbb P^1$.
\end{enumerate}
\end{thm}

\begin{thm}[Section \ref{kod=1}] \label{thm2}
Let $X$ be a projective  manifold whose canonical linear system
is composite with a pencil. Assume that
$$
q(X)>\op{min}\{\op{a}(X)+1, \dim(X)\},
$$
and that $\kappa(\op{alb}(X))=1$. Let $F$ denote the general fibre of
the canonical fibration $f: X\to C$.\\
Then $p_g(F)=1$, and $f$ is the composition of
$\op{alb}_X$ with the Ueno map $u_{\op{alb}(X)}$ (cf. (\ref{prop})).
Moreover $ \ p_g(X)+1 \geq g(C) \geq 2$, and $q(X)=g(C)+\op{a}(X)-1.$
\end{thm}

So the assumptions made in Theorem \ref{thm1}, (i), and in Theorem
\ref{thm2} imply that the image $C$ of the canonical fibration $f$ is a curve
positive genus, hence the universal property of the Albanese map implies that
the canonical fibration factors through the Albanese map. For threefolds
with $q(X)\geq 5$, one has a similar property.

\begin{thm}[Section \ref{gt}]\label{3thm}
Let $X$ be a threefold of general type, whose canonical linear system is composite
with a pencil. Assume either one of the following conditions:
\begin{enumerate}
\item[a.] $q(X)\geq 5$.
\item[b.] $q(X)\geq 4$ and $X$ is of maximal Albanese dimension.
\item[c.] $q(X)\geq 4$ and $\op{alb}(X)$ is a surface of Kodaira dimension $1$.
\item[d.] $q(X)\geq 3$ and $\op{alb}(X)$ is a curve.
\end{enumerate}
Then the Stein factorization of the Albanese map $X\to \op{alb}(X)$ factors
through the canonical fibration.
\end{thm}
Of course, parts b), c) and d) in \ref{3thm} are just \ref{thm1} and \ref{thm2},
respectively. We do not know, whether the image of the canonical fibration
can be $\BB P^1$ under the assumption a).

The bound $q(X)\geq 5$ is the best possible. In Section
\ref{gt} we give examples of projective $n$-folds $X$ of general
type with $q(X)=2n-2$ and
$$
\kappa(\op{alb}(X))= \op{a}(X)=n-1,
$$
whose canonical fibration is a morphism $f:X\to \BB{P}^1$. If $h:X\to Y$
denotes the Stein factorization of $\op{alb}_X$, then the induced map
$X\to \BB{P}^1\times Y$ is generically finite. In
those examples the geometric genus of a general fiber $F$ of $f$ is
$p_g(F)=2^{n-1}+1>1$.

It is an open question whether there exists a constant $c(n)$, depending
only $n$, such that for all $n$-dimensional manifolds with
$q(X) \geq c(n)$, $p_g(X)\geq 2$ and $a(X)< \dim(X)$, the morphism
$X \to Z \times \op{alb}(X)$ induced by the Stein factorization $X\to Z$
of the canonical map and the Albanese map, is never surjective.

In Section \ref{q=dim} we will consider projective manifolds of
general type and maximal Albanese dimension with $q(X)=\dim(X)$
whose canonical system is composite with a pencil (Examples
\ref{e1} and \ref{ex1}).

For a projective variety $X'$ we write $p_g(X')$, and $q(X')$ for
the geometric genus and the irregularity of any smooth model $X$
of $X'$, respectively. At several places we change the smooth
birational model $X$. In particular we will always assume that the
canonical map $\Phi_X$ is a morphism. We will say that the canonical fibration
is a (rational / irrational) pencil, whenever the canonical linear system
is composite with a pencil (and the image of the canonical fibration
a rational / non-rational curve).

\section{Manifolds of maximal Albanese dimension}\label{sect2}
Recall a well-known result on the geometry of subvarieties of an
Abelian variety.
\begin{prop}[Ueno \cite{Ue}, Theorem 10.9, and Koll\'ar, see \cite{Mo},
Corollary 3.10]\label{prop} Let $A$ be an Abelian variety, and
$Y\subset A$ an irreducible reduced subvariety. Let $A_0$ be the
connected component of $\{a\in A; \  a+ Y= Y\}$ containing the
origin. Then $A_0$ is an Abelian subvariety of $A$, and
$$u_Y: Y\>>>  B=Y/A_0 \subset A/A_0$$
is an \'etale fiber bundle with fiber $A_0$. Moreover $B$ is of
general type, $u_Y$ is the canonical map and it is
birational to the Iitaka fibration of $Y$.
\end{prop}
\noindent
We will call $u_Y: Y\to B=Y/A_0$ the Ueno map.

For the proof of Theorem \ref{thm1} we need two simple observations.
\begin{prop} \label{trace}
Let $\alpha: X \to Y $ be a generically finite surjective morphism of
projective manifolds. Then, $\omega_{Y}$ is a direct factor of
$\alpha_* \omega_{X}$.
\end{prop}
\begin{proof}
The statement is compatible with blowing up $X$ and $Y$, hence one
may assume the branch locus $\Delta$ of $\alpha$ to be a normal
crossing divisor. Let
$$
X\>{g}>> Z\>{h}>> Y
$$
be the Stein factorization of $\alpha$. Then $Z$ has at most quotient
singularities and hence rational singularities (cf. e.g., \cite{EV},
Lemma 3.24), and $g_*\omega_{X}=\omega_Z$.
Duality for finite morphisms (see \cite{Ha}, Ex. II 6.10,
and \cite{Ha}, Ex. III 7.2) implies
$$
\alpha_* \omega_{X}=h_*\omega_Z = \mathcal Hom_{Y}(h_*\mathcal O_Z, \omega_{Y}),
$$
so $\alpha_*\omega_X$ contains $\mathcal Hom_{Y}(\mathcal O_Y,
\omega_{Y})=\omega_Y$ as a direct factor.
\end{proof}

\begin{lem}\label{g2}
Let $X$  be a projective  manifold  with $p_g(X)\ge2$.
Assume that the  canonical fibration $f: X\to C$ is a pencil
with general fiber $F$.

If $g(C)\ge2$, then either $p_g(F)=1$ or $p_g(F)=p_g(X)=g(C)=2$.
\end{lem}
\begin{proof}
Since $p_g(X)\not=0$, we have $p_g(F)>0$. Assume that
$p_g(F)\neq 1$. Since $\Phi_X$ factors through $f$, the moving
part of $|K_X|$ is given by $f^*\Cal L$ for some invertible subsheaf
$\Cal L$ of $f_*\omega_X$, splitting locally and of strictly positive degree;
thus
$$p_g(X)=h^0(\Cal L)=h^0(f_*\omega_X)\text{ \ \ and \ \ }
h^0(\omega_C\otimes \Cal L^{-1}) < h^0(\omega_C).$$
Consider the exact sequence of sheaves
$$
0\>>>\Cal L\>>> f_*\omega_X\>>> \Cal Q\>>> 0.
$$
By the choice of $\Cal L$ the sheaf $\Cal Q$ is locally free
of rank $p_g(F)-1$. By \cite{Fu} the sheaf $f_*\omega_{X/C}=
f_*\omega_X\otimes \omega_C^{-1}$ is nef, hence $\deg(\Cal Q \otimes
\omega_C^{-1})\geq 0$. The Riemann-Roch Theorem for locally free sheaves
on $C$ implies that
\begin{gather*}  g(C)-1\ge h^0(\omega_C\otimes \Cal L^{-1}) =
h^1(\Cal L) \ge h^0(\Cal Q) \cr \ge (p_g(F)-1)(g(C)-1).
\end{gather*}
Then $p_g(F)=2$ and $h^1(\Cal L)=g(C)-1$, hence $h^0(\Cal L)=\deg(\Cal L)$.
By Clifford's Theorem $\deg (\Cal L)=2$ and
$\Cal L\simeq \omega_C$. So $p_g(X)=g(C)=2$.
\end{proof}

\begin{proof}[Proof of Theorem \ref{thm1}]
Blowing up $X$, if necessary, we may assume that
the Albanese map  $\op{alb}: X\to \op{alb}(X)
\subset\op{Alb}(X)$ factors like
$$
X \> \alpha >> Y \> \delta >> \op{alb}(X),$$ for some
desingularization $Y$ of $\op{alb}(X)$. In particular, by
\ref{trace} there exists some torsion-free sheaf $\mathcal G$ with
$\alpha_*\omega_X = \omega_Y \oplus \mathcal G$. So the canonical
map $\Phi_Y$ factors through $\Phi_X$, hence through the canonical
fibration $f:X \to C$, and the image of $\Phi_Y$ must be a curve
$B$. By \ref{prop} $\Phi_Y: Y \to B$ is birational to an \'etale fibre
bundle of Abelian subvarieties of $\op{Alb}(X)$ of dimension
$\dim(Y)-1=\dim(X)-1$, and $q(X)=g(B)+\dim(X)-1$. On the other
hand, the general fibre $F$ of $f$ is of maximal Albanese
dimension, and its image in $\op{Alb}(X)$ is an Abelian
subvariety. Hence $q(X) \geq g(C) + \dim(F)$ and $g(C)\leq g(B)$.
Since $C$ is a covering of $B$, this is only possible
for $B=C$ and $g(C)=g(B)\geq 2$. Let us write $\iota: Y \to C$, such
that $f=\iota\circ \alpha$. By Lemma \ref{prop}, if $p_g(F) \neq 1$ one has
$$p_g(F)=p_g(X)=g(C)=2\text{ \ \ and \ \ }q(X)=\dim(X)+1.$$
So $\op{rank}f_*\omega_{X}=2$, and $\iota_*\mathcal G$ is an
invertible sheaf. Since $f_*\omega_{X/C}$ is nef, $\iota_*\mathcal
G\otimes \omega_C^{-1}$ is nef, and the Riemann-Roch Theorem
implies $h^0(i_*\mathcal G)>0$, or
$$p_g(X)=h^0(i_*\omega_Y)+h^0( i_*\mathcal G)\geq3,$$
a contradiction.

Since $F$ is of maximal Albanese dimension $p_g(F)=1$ implies that $q(F)=\dim(F)$
and $\chi(\mc O_F)=0$ (see \cite{Ue}). The restriction map
$$r: \op{Pic}^0(X)\>>>  \op{Pic}^0(F)$$
is surjective. By the generic vanishing theorem (\cite{GL} or \cite{EL}),
for a general $\mc L \in \op{Pic}^0(X)$ the sheaf $\omega_F\otimes\mc L|_F$
has no higher cohomology, hence
$$h^0(\omega_F\otimes\mc L|_F)= \chi(\omega_F\otimes\mc L|_F)=0.$$
Tensoring the exact sequences
$$0\>>>  \omega_X(-(n+1)F)\>>>  \omega_X(-nF)\>>> \omega_F\>>> 0$$
with $\mc L$, one obtains by descending induction on $n$ that
$h^0(\omega_X\otimes\mc L)=0$. Using again the generic vanishing theorem
(\cite{GL} or \cite{EL}) one finds $\chi(\mc O_X)=0$.
\end{proof}

When $\dim(X)=3$, hence $\dim(F)=2$, the condition $\chi(\mc O_F)=0$
forces the Kodaira dimension of $F$ to be strictly smaller than two.
By the easy addition formula for the Kodaira dimension one obtains an
improvement of \cite{Ca}, Theorem 1.
\begin{cor}\label{impr}
Let $X$ be a projective  $3$-fold of general type and maximal
Albanese dimension. If  $q(X)\geq 4$, then  the  canonical linear
system of $X$ is not composite with a pencil.
In particular, if $p_g(X)=2$ then $q(X)=3$.
\end{cor}

\begin{cor}\label{cor2}
Let $X$ be a projective manifold of maximal  Albanese dimension.
If $p_g(X)=2$, then $q(X)\leq \dim(X)+1$.
\end{cor}
\begin{proof}
Let  $f: X\to C$ be the canonical fibration.
If $q(X)>\dim(X)+1$, Theorem \ref{thm1} implies that
$p_g(X)= g(C)=q(X)-\dim(X)+1\geq 3$.
\end{proof}

\begin{rk}
In dimension greater than three there are manifolds $X$ satisfying
the assumptions made in Theorem \ref{thm1} or Corollary
\ref{cor2}. For example, by \cite{CH2} there exist smooth
manifolds $F$ of general type and of maximal Albanese dimension
with $p_g(F)=1$. If $C$ is a curve of genus $\geq 2$ then one can
choose $X=F\times C$.
\end{rk}

The arguments used to prove Lemma \ref{g2} can also be applied
to bicanonical maps.
\begin{thm} Let $X$ be a projective  manifold of maximal  Albanese
dimension, with $\kappa (X)\geq 2$ and $q(X)>\dim(X)$.
Then the bicanonical linear system is not composite with a pencil.
\end{thm}
\begin{proof} Suppose that $\Phi_{2X}$ is composite with a pencil,
and let $h: X\to B$ be the fibration obtained as the Stein
factorization of $\Phi_{2X}$. A general fiber $F$ of $h$ has
maximal Albanese dimension, and by the weak addition theorem
$\kappa(F) \geq 1 $, in particular $F$ is not an Abelian variety.
By \cite{CH1}, Theorem 3.2, one finds $p_2(F)\geq 2$.

As in the proof of Lemma \ref{g2}, there is an invertible subsheaf
$\Cal L \hookrightarrow h_*\omega_X^{\otimes2}$, splitting locally
and with $h^0(\Cal L)=p_2(X)$. The locally free quotient
sheaf $\Cal Q = h_*\omega_X^{\otimes2}/\Cal L$ is locally free
of rank $p_2(F)-1$ and
$$g(B)-1 \ge  h^0(\omega_B\otimes \Cal L^{-1}) =h^1(\Cal L) \geq h^0(\Cal Q).$$
By \cite{V} the sheaf $h_*\omega_{X/B}^{\otimes2}$ is nef, hence
$$\deg(\Cal Q)-2(p_2(F)-1)(2g(B)-2)= \deg(\Cal Q\otimes
\omega_B^{-2})\geq 0.$$
The Riemann-Roch Theorem implies
$$h^0(\Cal Q)\geq 3(p_2(F)-1)(g(B)-1),$$
and hence $g(B)\leq 1$.

By \ref{prop} and \ref{trace} one has $p_g(X) \geq 2$.
The canonical map $\Phi_X $ is defined by a linear subsystem of
$|K_{2X}|$, hence it factors through $\Phi_{2X} $,
and the canonical fibration $f:X \to C$ factors through $h$.
Since the fibres of $f$ and $h$ are both connected of the same dimension,
$C=B$ is a curve of genus $\leq 1$, contradicting Theorem \ref{thm1}.
\end{proof}

\section{The case $\kappa(\op{alb}(X))=1$ }\label{kod=1}

\begin{lem}\label{prod}
Let $X$ and $B$ be projective manifolds. Assume that the canonical
fibration is a pencil $f:X\to C$, and let $h:X\to B$ be a morphism, such that
$\alpha=(f,h):X\to C\times B$ is surjective. Then there is an injection
$$\mc O_B^{\oplus p_g(X)} \>>> h_*\omega_X.$$
\end{lem}
\begin{proof}
Let us write $\alpha=(f, h): X\to C\times B$.
Since $\Phi_X$ factors through $f$, we may choose again an invertible
subsheaf $\mc L$ of $f_*\omega_X=pr_{1*}\alpha_*\omega_X$
with $p_g(X)=h^0(\mc L)$. Then $pr_1^*\mc L$ is a subsheaf of
$\alpha_*\omega_X$ and one finds
$$
\mc O_B^{\oplus h^0(\mc L)} \simeq pr_{2*}pr_1^*\mc L \> \subset >>
pr_{2*}\alpha_* \omega_X=h_*\omega_X.
$$
\end{proof}
\begin{cor}\label{nonex}
Let  $X$ be a projective  manifold, whose canonical linear system
is composite with a pencil, and let $h: X\to B$ be a fiber space
over a curve $B$ with $g(B)\geq 3$. Then $h$ factors through the canonical
fibration $f:X\to C$.
\end{cor}
\begin{proof}
If not, the induced map $\alpha=(f,h):X \to C\times B$ is surjective
and \ref{prod} implies that $\op{rank}(h_*\omega_X)\geq p_g(X)$.
Using again the nefness of $h_*\omega_{X/B}$, the Riemann-Roch
theorem implies that
\begin{gather*}
p_g(X)=h^0(h_*\omega_X) \geq \deg(h_*\omega_{X/B})+
\op{rank}(h_*\omega_X)(g(B)-1)\\
\geq p_g(X)(g(B)-1),
\end{gather*}
for $g(B)\geq 3$ a contradiction.
\end{proof}

\begin{proof}[Proof of Theorem \ref{thm2}]
By Theorem \ref{thm1}, we can assume that $\op{a}(X)<\dim(X)$.
Consider
$$\begin{CD} h=u\circ\op{alb}:  X \>\op{alb}>> \op{alb}(X)
\>u>> B,\end{CD} $$ where $u$ is the Ueno map. Then $B$ is a
curve of genus
$$
g(B)=q(X)-\op{a}(X)+1\geq 3.
$$
By Corollary \ref{nonex} $h$ factors through the canonical fibration $f:X\to C$.
Then the image of a general fibre $F$ of $f$ under $\op{alb}$ is an Abelian
variety of dimension $a(X)-1$, hence $ g(C) \leq q(X) -a(X) + 1$,
and $C=B$. By Lemma \ref{g2}, one has $p_g(F)=1$. Since
$f_*\omega_{X/C}$ is invertible and nef, the Riemann-Roch theorem
implies that
$$
p_g(X)=h^0(f_*\omega_X) \geq \deg(f_*\omega_{X/C})+
g(C)-1\geq g(C)-1.$$
\end{proof}

\begin{rk}
Examples, due to Beauville (\cite{Be}, Example 2), show that the
condition on $q(X)$ in Theorem \ref{thm2} can not be weakened.
\end{rk}
There is a series of examples, showing that there are no further
restrictions on the numerical invariants in \ref{thm2}.
\begin{exa} Let $T$ be a surface of general type
with $p_g(T)=q(T)=0$ and $\pi_1(T)=<\sigma>\simeq \mathbb Z_2$.
Let $S\to T$ be the universal cover. Consider in addition a curve
$\ti C$ with an involution $\tau$ such that the genus of $C=\ti
C/\tau$ is $\geq 3$, and let $R$ be the ramification divisor of $\ti C$
over $C$. Then
$$
X=(S\times \ti C)/\sigma\times \tau
$$
is of general type with $p_g(X)=g(C)+\frac{1}{2}\deg(R)-1$,
$q(X)=g(C)$, and $$\Phi_X=\op{alb}_X: X \>>>  C=\ti C/\tau$$ is
induced by the second projection.
\end{exa}
As we will see in Section \ref{gt}, for $\dim(X)=3$ and
$q(X)=4$ the condition on the Kodaira dimension of $\op{alb}(X)$
in Theorem \ref{thm2} is needed. Nevertheless one can show
without this assumption:
\begin{prop}
Let $X$ be a projective  $3$-fold of general type, whose canonical
fibration $f: X\to C$ is a pencil, and let $F$ be a general fiber of $f$. If
$q(X)\geq 4$ and $g(C)\geq2$, then $p_g(F)=1$, and $f$ is the
composite of $\op{alb}_X$ and the Ueno map $u_{\op{alb}(X)}$.
\end{prop}
\begin{proof}
By Theorems \ref{thm1} and \ref{thm2}, we may assume that
a desingularization $S$ of $\op{alb}(X)$ is a surface
of general type. By Lemma \ref{g2}, it only
remains to exclude the case $g(C)=p_g(X)=p_g(F)=2$.

By the universal property of the Albanese map, there exists
a morphism $h: S\to C$. Since $f=h\circ\op{alb_X}$,
the fibres of $h$ are connected.

For a general point $p\in C$, consider the fibres $F=f^*(p)$,
and $H=h^*(p)$. Since $X$ and $S$ are of general type, the easy
addition formula implies that $F$ and $H$ are of general type. So
$$ 2=p_g(F)\geq q(F)\geq g(H)\geq2,$$
and $q(F)=g(H)=2$. Then
$$4\leq q(X)\leq q(F)+g(C)=4 \text{ \ \ and \ \ }
q(S)=q(X)=4=g(H)+g(C).$$
This implies that $S\simeq C\times H$. Let $\alpha: X\to H$ be the
composite of $\op{alb_X}$ with the projection $S\to H$.

The sheaf $\mc E=\alpha_*\omega_{X/H}$ is nef.
Assume that for $p\in H$ in general position
$$H^0(\mc E\otimes\omega_H\otimes \mc O_H(-p) )\neq 0.$$
Then $\alpha^{-1}(p)$ belongs to the moving part of the linear system
$|K_X|$, hence it must be a fibre of $f$, contradicting
the construction of $\alpha$.

Using the Riemann-Roch theorem one finds
\begin{gather*}
0=h^0(\mc E\otimes\omega_H\otimes \mc O_H(-p))=
\deg(\mc E) + h^1(\mc E\otimes\omega_H\otimes \mc O_H(-p)),
\end{gather*}
and both, $\deg(\mc E)$ and $h^1(\mc E\otimes\omega_H\otimes \mc O_H(-p) )$
are zero.
The long exact cohomology sequence for
$$0\>>> \mc E\otimes\omega_H\otimes \mc O_H(-p)\>>>
\mc E\otimes\omega_H\>>> \mathbb C_p^{\op{rank}(\mc E)}
\>>> 0,$$
implies that $\op{rank}(\mc E)=p_g(X)=2$, hence
$p_g(\alpha^{-1}(p))=2$, and
$\Phi_{X|\alpha^{-1}(p)}$ is the canonical map of $\alpha^{-1}(p)$.
Then the image of the canonical fibration of $\alpha^{-1}(p)$
is a curve of genus $\geq 2$, contradicting \cite{X1}.
\end{proof}

\section{$\op{alb} (X)$ is of general type }\label{gt}

As mentioned in the Introduction, Theorem \ref{thm2} does not hold
true, without the assumption $\kappa(\op{alb} (X))=1$.

\begin{exa}\label{q=2d-2} For $n\geq3$ and
$$\frac{3n-3}{2} \leq q \leq 2n-2,$$
we will construct projective $n$-folds $X$  with irregularity $q(X)=q$
and with $\kappa(\op{alb}(X))= \op{a}(X)=n-1$,
whose canonical linear system is composite with a rational pencil with
connected fibres. To this aim, consider an \'etale double cover
$\ti Y\to Y=\ti Y/\sigma$, satisfying the following condition:
\begin{enumerate}
\item[($*$)] $Y$ is a manifold of general type, of maximal Albanese
dimension, and with $p_g(\ti Y)=p_g(Y)+1$.
\end{enumerate}
Let $C$ be a smooth curve of genus $g\geq 2$ with hyperelliptic
involution $\tau$. Choose
$$
X=(\ti Y\times C)/\sigma\times \tau,
$$
and $f: X\to \mathbb P^1=C/\tau$ to be the fibration induced
by the second projection. Then $X$ is of general type with
$$p_g(X)=g, \ \ q(X)=q(Y), \ \ \op{alb}_X(X)=\op{alb}_Y(Y),$$
and with $\Phi_X=f:X\to \BB P^1$. The product of the canonical fibration
and the Albanese map factors like
$$
X \>>> (\ti Y/\sigma)\times (C/\tau) = \BB P^1\times Y \>>> \BB P^1 \times
\op{alb}_Y(Y),
$$
hence it is surjective.\\

To show that there exist examples of $\ti Y\to Y$
satisfying the condition ($*$) and with $q(Y)=q$ we start with:
\begin{enumerate}
\item[(i)] $\dim(Y)=1$: Consider an \'etale cover $\ti Y\to Y$ of
a curve $Y$ of genus $2$, hence $p_g(Y)=2$, $p_g(\ti Y)=3$, and
$q(Y)=2$.
\item[(ii)] $\dim(Y)=2$: Let $B$ be a non-hyperelliptic
curve of genus $3$, let $Y=B^{(2)}\subset \op{Jac}(B)$ be the
Theta divisor, and let $\alpha : \bar A\to \op{Jac}(B)$ be the
double cover defined by a non-zero $2$-torsion element of
$\op{Pic}^0\op{Jac}(B)$. Then the inverse image $\ti Y\subset A$
of $Y$ under $\alpha$ is an \'etale double cover of $Y$, with
$p_g(Y)=3$, $p_g(\ti Y)=4$ and $q(Y)=3$.
\item[(iii)] Assume that
$\ti Y_i\to Y_i$ ($i=1$, $2$) satisfies ($*$), and let $\sigma_i$
be the corresponding involution. Then
$$\ti Y:=(\ti Y_1\times
\ti Y_2)/(\sigma_1\times\sigma_2) \>>>  Y_1\times Y_2$$ is an
\'etale double cover with $p_g(\ti Y)=p_g(Y_1)p_g(Y_2)+1$ and
$q(\ti Y)=q(Y_1)+q(Y_2)$, hence it also satisfies ($*$).
\end{enumerate}
Now given $n\geq 3$ taking for $Y$ $n-1$ copies of (i) we obtain examples
with $q(X)=2n-2$. If $n$ is odd, we may choose $\frac{n-1}{2}$ copies
of (ii) and we find $X$ with $q(X)=\frac{3n-3}{2}$. For $n$ even,
one reaches $q(X)=\frac{3n-2}{2}$. Of course, one can also find examples for all
intermediate values.
\end{exa}

\begin{proof}[Proof of Theorem \ref{3thm}]
Remark that part b) is a special case of Theorem \ref{thm1}
and that c) and d) are special cases of Theorem \ref{thm2}.
Moreover, if the canonical fibration $f:X\to C$ is an irrational pencil,
the Albanese map factors through $f$.

Hence it only remains to consider the case
$$
C=\BB P^1, \ \ \ \
q(X) \geq 5, \text{ \ \ and \ \ }\
\dim(\op{alb}(X))=\kappa(\op{alb}(X)) = 2.
$$
Let
$$
X\>{h}>> S \>>> \op{alb}(X)
$$
be the Stein factorization of $\op{alb}_X$. Blowing up $X$, if necessary,
we will assume that $S$ is non-singular. Remark that $q(S)=q(X)$.

We will show step by step, that the surjectivity of
$$
\alpha=(f,h):X\>>> \BB P^1\times S
$$
implies:
\begin{enumerate}
\item[(i)] The evaluation map $H^0(h_* \omega_X)\otimes \mc O_S
\to h_*\omega_X$ is injective; in particular,
$\op{rank}(h_*\omega_X)\geq p_g(X)$.
\item[(ii)] $h^0(\Omega_X^2)=p_g(S)$.
\item[(iii)] $h^i(h_*\omega_X)=0$ for $i>0$.
\item[(iv)] $\chi(\mc O_S)=0$.
\end{enumerate}
Of course, (iv) contradicts the assumption that $\op{alb}(X)$, hence $S$,
is a surface of general type.

(i) follows from Lemma \ref{prod}.

For (ii), replacing $S$ and $X$ by suitable smooth
birational models, we will assume that for some simple
normal crossing divisor $B$ on $S$ the morphism
$h$ is smooth over $S\setminus B$, and that
$\Delta =(h^*B)_{red}$ is a divisor with simple normal crossings.
We have morphisms of sheaves
\begin{gather*}
h^*\Omega_S^1\> \subset >> \mc E=\Omega_X^1\cap
h^*\Omega_S^1(\op{log}B) \>\subset >>  \Omega_X^1
\>{\epsilon}>>  \omega_{X/S}(\Delta-h^*B),
\end{gather*}
where $\epsilon$ is the composite
$\Omega_X^1\to\Omega_X^1(\op{log}\Delta)\to
\omega_{X/S}(\Delta-h^*B) $.
Taking wedge product, one finds a complex of sheaves
$$0 \>>>  h^*\omega_S\>>> \Omega_X^2\>>>  \omega_{X/S}
\otimes h^*\Omega_S^1.$$
Note that $h_*(\mc O_X(h^*B-\Delta))=\mc O_S$, and that
$h^0(h^*\omega_S)=h^0(\op{det}(\mc E))$.

Now suppose that $h^0(\Omega_X^2)>p_g(S)$,
choose $\varphi\in H^0(\Omega_X^2)\setminus h^*H^0(\omega_S) $,
and let
\begin{equation}\label{1}
\mc O_X\>>>  \omega_{X/S}
\otimes h^*\Omega_S^1.
\end{equation}
be the induced non-trivial map.
Tensoring (\ref{1}) with $h^*\Omega_S^1$, one obtains
$$h^*\Omega_S^1\>>>  \omega_{X/S}
\otimes h^*\Omega_S^1\otimes h^*\Omega_S^1\>>>
\omega_{X/S}\otimes h^*\omega_S$$
and, taking the direct image, one has
a non-trivial map
$$r: \Omega_S^1\>>>  h_*\omega_X.$$
Let $\mc K=\op{Ker} (r : H^0(\Omega_S^1)\to H^0( h_*\omega_X)).$
Clearly $$\mc K=\{s\in H^0(\Omega_S^1)| \ s\wedge \varphi=0\}.$$
Since $\op{rank} (\op{Im}\ r)\leq2$,
(i) implies that $h^0(\op{Im}\ r)\leq2$ and
$\dim(\mc K)\geq q(S)-2$.

For a general point $x \in X$ define
$$
\Sigma=\{s\in\Omega_{x,X}^1| \ s\wedge \varphi=0 \}.
$$
Clearly $\dim(\Sigma)\leq2$,
and by the choice of $\varphi$,
$$\Sigma\not \subset h^*\Omega_{h(x),S}^1.$$
This implies that $h^*\mc K$ generates
an invertible subsheaf of $h^*\Omega_S^1$, hence a morphism
$S\to D$ to a curve $D$ of genus $g(D)=\dim(\mc K)
\geq q(S)-2\geq3$, contradicting \ref{nonex}.

The Leray spectral sequence and the Serre duality for $h$ imply that
$$h^i(\omega_X)=h^i(h_*\omega_X) + h^{i-1}(R^1h_*\omega_X)=
h^i(h_*\omega_X) + h^{i-1}(\omega_S).
$$
Hence (iii) follows from $q(X)=q(S)$ and from (ii).

For (iv) remark first that (iii) together with semicontinuity, imply
that for some open neighborhood $U \subset \op{Pic}^0(S)$ of $\mc O_S$,
and all $\mc P \in U$ one has
$$
H^i(h_* \omega_X \otimes \mc P)=0\text{ \ \ for \ \ } i>0.
$$
It follows that for $\mc P \in U$,
$$
h^0(h_*\omega_X \otimes \mc P)=\chi(h_*\omega_X)=p_g(X).$$

Let $\mc Q$ be the cokernel of the injection
$\mc O_S^{\oplus p_g} \to h_*\omega_X $
constructed in (i). For $\mc \mc P \in U$ one finds an exact sequence
$$ 0 \to H^0(\mc O_S^{\oplus p_g} \otimes \mc P) \to H^0(h_*\omega_X  \otimes \mc P)
\to H^0( \mc Q  \otimes \mc P) \to H^1(\mc O_S^{\oplus p_g}  \otimes \mc P) \to 0,$$
an isomorphism
$$
H^1( \mc Q  \otimes \mc P) \>>>  H^2(\mc O_S^{\oplus p_g}  \otimes \mc P)
$$
and $H^2( \mc Q  \otimes \mc P)=0$. For $\mc P=\mc O_S$, one finds
\begin{gather*}
h^0(\mc Q)= h^1(\mc O_S^{\oplus p_g})=p_g(X) q(X),\\
h^1(\mc Q)= h^2(\mc O_S^{\oplus p_g})=p_g(X) p_g(S),\\
h^2(\mc Q)=0,
\end{gather*}
whereas for $\mc P \ne \mc O_S$,
\begin{gather*}
h^0(\mc Q  \otimes \mc P)=h^0(h_*\omega_X  \otimes \mc P)= h^0(h_*\omega_X
\otimes \mc P)=p_g(X),\\
h^1(\mc Q  \otimes \mc P)=h^2(\mc Q  \otimes \mc P)=0.
\end{gather*}
Since $\chi(\mc Q)=\chi(\mc Q \otimes \mc P)$, one obtains for the
Euler characteristic
$$
\chi(\mc Q)=p_g(X) q-p_g(X)p_g(S)=p_g(X).
$$
Therefore, $q(S)=p_g(S)+1$, and hence $\chi(\mc O_S)=0$.
\end{proof}
\begin{rk}\label{jungkai}
In a first version of this note, we were not able to show that
(i)-(iii) implies (iv), hence not able to prove Theorem \ref{3thm}, a), as
it is stated now. The argument presented above, has been told us by
Jungkai Alfred Chen.
\end{rk}
\section{The case $a(X)=q(X)=\dim(X)$ }\label{q=dim}
\begin{prop}\label{sur} Let $X$ be a projective  manifold
whose canonical system is composite with a pencil, and
let $f: X\to C$ be the canonical fibration.
If the Albanese map $\op{alb}: X\to \op{Alb}(X)$
is surjective, then either $C\simeq \mathbb P^1$,
or $g(C)=1$ and $q(X)<\dim(X)$.
\end{prop}
\begin{proof}
Clearly $g(C)\leq1$. If $C$ is an elliptic curve,
consider the commutative diagram
$$\begin{CD} X \>\op{alb}>> A:=\op{Alb(X)} \\
\V f VV \V h VV    \\
C\>=>> C
\end{CD} $$
where $h$ is induced by the universal property of Albanese map.
If $q(X)=\dim(X)$, Proposition \ref{trace} implies that
$\op{alb}_*\omega_X=\omega_A\oplus \mc F$
for some torsion-free sheaf $\mc F$ on $A$. Let $\mc E \subset
f_*\omega_X=h_*\omega_A\oplus h_*\mc F$ be the subsheaf generated
by the global sections of $f_*\omega_X$. Since $p_g(X)>p_g(A)=1$,
one finds that $h^0(\mc F)>0$, hence $\op{rank}(\mc E) \geq 2$,
contradicting the assumption that $\Phi_X$ factors through $f$.
\end{proof}

\begin{cor}\label{cor} Let $X$ be a projective manifold,
with $q(X)=\dim(X)$ and of maximal  Albanese dimension.
Assume the canonical system is composite with a pencil, and
let $f: X\to C$ be the canonical fibration.
Then $C\simeq \mathbb P^1$.
\end{cor}

Let us end by giving some examples of manifolds $X$
of maximal Albanese dimension with $q(X)=\dim(X),$
and $p_g(X)$ arbitrarily large,
whose canonical map is composite with a rational pencil.

\begin{exa}\label{e1}
Let ($A, \Theta$) be a principally polarized Abelian $n$-fold
($n\ge3$), and $D\in |2\Theta|$  a smooth divisor. The double
cover obtained by taking the root out of $D$ is of general type with
$p_g(X)=2$.
\end{exa}

\begin{exa}[For $n=2$ see \cite{X2}]\label{ex1}
For $n \ge 2$ let ($A, \Theta$) be a simple polarized Abelian $n$-fold
of type  $(1, \cdots , 1, 2)$. Let $\sigma : X'\to A$ be a
blowing up, such that the moving part of
$|\sigma^*\Theta|$ defines a morphism $f': X'\to C'$, necessarily with
$C'=\mathbb P^1$. The general fibre $F'$ of $f'$ is of general type
and $p_g(F')=h^0(\omega_{\Theta})=n+1$. Since $p_g(X')=1$, one may
write
\begin{equation}\label{2}
f'_*\omega_{X'}=\mc O_{\BB P^1}\oplus \mc O_{\mathbb P^1}(a_1) \oplus \cdots
\oplus \mc O_{\mathbb P^1}(a_n),
\end{equation}
with $-1 \geq a_1 \geq a_2 \geq \ldots \geq a_n$. Since
$f'_*\omega_{X'/\mathbb P^1}$ is nef by \cite{Fu}, $a_n\geq -2$.
On the other hand, $\sigma_*(\omega_{X'}(F'))$ is a subsheaf of
$\mc O_A(\Theta)$, hence $h^0(\omega_{X'}(F'))\leq 2$.
One finds in (\ref{2}) $a_1=a_2=\cdots =a_n =-2$.

Let $\pi :  C=\mathbb P^1\to C'$ be a finite morphism of degree
$d>1$, ramified only over points $p$ with ${f'}^{-1}(p)$ smooth,
and let $f:  X\to  C$ be the pull-back family.

Base change and (\ref{2}) imply that
$$
f_*\omega_{ X/C}=\pi^*f'_*\omega_{X'/ C'}= \mc O_{\mathbb P^1}(2d)
\oplus \mc O_{\mathbb P^1}^{\oplus n},
$$
hence that $f_*\omega_{ X}= \mc O_{\mathbb P^1}(2d-2)
\oplus \mc O_{\mathbb P^1}(-2)^{\oplus n}.$
Then $p_g( X)=2d-1$, the canonical map $\Phi_{ X}$ factors through $f$, and
$\omega_{X}$ contains $f^*\mc O_{\mathbb P^1}(2)$. The latter induces
an inclusion $\omega_X/C\to \omega_X^2$, and since
the general fiber of $f$ is of general type, $X$ must be of general type.

By construction $X$ is of maximal Albanese dimension, and Theorem \ref{thm1}
implies that $\dim(X)=q(X)=q(X')$. This of course follows by base change
for $f'_*\mc O_{X'}$ and $R^1f'_*\mc O_{X'}$, as well.
\end{exa}

{\it Acknowledgments.} Part of this work was done during the first
author's stay at the Universit\"at Essen. He thanks the members of
the Department of Mathematics, in particular H\'el\`ene Esnault,
for their help and hospitality. We are both very grateful to
Jungkai Alfred Chen who realized that the vanishing of the higher cohomology
of the sheaves $h_*\omega_X$ in the proof of Theorem \ref{3thm}, a), can not
hold true if $S$ is a surface of general type, and who allowed us to
include his argument in this note.

\end{document}